\documentstyle[amstex,12pt]{article}
\newcommand{\lon}{\longrightarrow}
\newcommand{\rar}{\rightarrow}
\newcommand{\hook}{\hookrightarrow}
\newcommand{\Proof}{\noindent{\bf Proof}.\, }

\newcommand{\R}{{\Bbb R}}

\newcommand{\p}{{\partial}}

\newcommand{\K}{{\Bbb K}}
\newcommand{\C}{{\Bbb C}}

\newcommand{\Beq}{\begin{equation}}
\newcommand{\Eeq}{\end{equation}}
\newcommand{\Beqr}{\begin{eqnarray*}}
\newcommand{\Eeqr}{\end{eqnarray*}}
\newcommand{\Beqrn}{\begin{eqnarray*}}
\newcommand{\Eeqrn}{\end{eqnarray*}}
\newcommand{\Ba}{\begin{array}}
\newcommand{\Ea}{\end{array}}
\newcommand{\Barr}{\begin{array}}
\newcommand{\Earr}{\end{array}}
\newcommand{\Bi}{\begin{itemize}}
\newcommand{\Ei}{\end{itemize}}
\newcommand{\Bc}{\begin{center}}
\newcommand{\Ec}{\end{center}}
\newcommand{\f}{{\cal O}}

\newcommand{\cJ}{{\cal J}}

\newcommand{\cQ}{{\cal Q}}
\newcommand{\cP}{{\cal P}}
\newcommand{\cN}{{\cal N}}

\newcommand{\al}{\alpha}

\newcommand{\la}{\lambda}

\newcommand{\vst}{\vspace{2 mm}}

\newcommand{\vse}{\vspace{8 mm}}
\begin{document}

\title{Deformation quantization of the $n$-tuple point}

\author{ S.A.\ Merkulov}
\date{}
\maketitle

\begin{abstract}
Contrary to the classical methods of quantum mechanics, the deformation quantization
can be carried out  on phase spaces which are not even topological manifolds.
In particular, the Moyal star product  gives rise to a canonical functor $F$ from
the category of  affine analytic spaces
to the category of associative
(in general, non-commutative) $\C$-algebras. Curiously, if $X$ is the
$n$-tuple point, $x^{n}=0$, then $F(X)$ is the algebra  of $n\times n$ matrices.
\end{abstract}

\vse

\paragraph{\bf 1. Introduction.}
This short note, which is largely about an entertaining
interpretation of the classical algebra of $n\times n$-matrices as a
quantized $n$-tuple point, is almost a mathematical anecdote.
This is also an attempt to understand what a quantum mechanical
system may be on spaces like the ``cross" $X_1=\{(x,y)\in \R^2 \mid
xy=0\}$, the ``tick" $X_2=\{(x,y)\in \R^2 \mid y^2-x^3=0\}$ or the
real line with one double point $X_3=\{(x,y)\in \R^2 \mid xy=0, y^{2}=0\}$ which
either fail to be topological manifolds or/and have nilpotents in their
structure sheaves. Contrary to the standard methods of quantum mechanics,
the deformation quantization \cite{F} (see also \cite{S} for an up-to-date
overview) easily sustains the
introduction of this type of singularities and equippes the (complexified) structure
sheaves of the associated
phase spaces with well-defined one-parameter non-commutative associative
star products $*_{\hbar}$ which, however, depend {\em
meromorphically}\, on the Planck constant $\hbar$. Their physical
interpretation is left to the imagination of the reader.

\vst

%%%%%%%%%%%%%%%%%%%%%%%%%%%%%%%%%%%%%%%%%%%%%%%%%%%%%%%%
\paragraph{2. The Moyal product $*_{\hbar}$ on affine
spaces.} Let $X$ be a subspace of $\R^n$ given by the equations
$$
\phi_{\al}(x)=0, \  \ \al=1,\ldots,k,
$$
where $\phi_{\al}(x)$ are polynomial (or analytic or even smooth)
functions on $\R^n$. The natural
coordinates on $\R^n$ are denoted by $x^a$, $a=1,\ldots,n$.
We understand the affine space $X$ as
a ringed space, i.e.\ as a pair $(X,\f_X)$ consisting of the subset of
points (with the induced topology),
 $X\subset \R^n$, satisfying the
above equations together with the structure sheaf $\f_X=\f_{\R^{n}}/J_X$,
where $J_X$ is the ideal subsheaf of the sheaf $\f_{\R^n}$ of  smooth
functions on $\R^n$ generated by $\phi_{\al}(x)$. Remarkably, the deformation
quantization on $X$ will {\em not}\, depend on the particular choice
of generators $\phi_{\al}$ of $J_X$ giving therefore rise to a
genuine functor on the category of affine spaces.

\vst

Let $M=\R^{2n}$ be the total space of the cotangent bundle to $\R^n$
with its canonical symplectic form $\omega=\sum_{a=1}^n dp_a\wedge dx^a$,
where $p_a$ are the natural fibre coordinates. The Moyal star product
\cite{M,F} makes the sheaf $\f_M$ of smooth functions on $M$,
or more precisely its extension $\f_M[[\la]]$, $\la$ being the formal
deformation parameter, into the sheaf of  non-commutative associative
algebras with the product given by
$$
f*_{\la} g := \left. e^{\sum_{a=1}^n \la\left(\frac{\p^2}{\p p_a\p
\tilde{x}^a} - \frac{\p^2}{\p x^a\p \tilde{p}_a}\right)}
f(x^b,p_b) g(\tilde{x}^c,\tilde{p}_c)\right|_{{x^a=\tilde{x}^a}\atop
p_a=\tilde{p}_a}.
$$
In the context of quantum mechanics the parameter $\la$ is set to be
$\frac{i}{2}\hbar$, $\hbar$ being the Planck constant, and the Moyal
product is denoted by $*_{\hbar}$.

\vst

With the affine subspace $(X, \f_X)$ of $\R^n$ we associate
\begin{itemize}
\item two subsheaves of ideals of $(\f_M[[\la]], *_{\la})$, the right ideal
$$
\cJ_r := \pi^*(J_X)*_{\la} \f_{M}[[\la]]
$$
and the left one
$$
\cJ_l := \f_M[[\la]]*_{\la} \pi^*(J_X),
$$
where $\pi:M=\Omega^1\R^n \rar \R^n$ is the natural projection;
\item two subsheaves of normalizers,
$$
\cN_r:= \{f\in \f_M[[\la]] \mid  f*_{\la} \pi^*(J_X) \subset \cJ_r\}
$$
and
$$
\cN_l:=\{ f\in \f_M[[\la]] \mid \pi^*(J_X)*_{\la} f \subset \cJ_l\},
$$
which are subsheaves of subrings of $(\f_M, *_{\la})$;
\item and, since $J_r\subset (\cN_r, *_{\la})$ and $J_l\subset (\cN_l,
*_{\la})$ are subsheaves of two-sided ideals, the two quotient sheaves
of (in general, non-commutative) associative algebras
$$
\left(\cP_X=\cN_r/\cJ_r, *_{\la}\right) \ \ \mbox{and} \ \
\left(\cQ_X=\cN_l/\cJ_l, *_{\la}\right).
$$
\end{itemize}

The star products in the sheaves $\cP_X$ and $\cQ_X$ are naturally
induced from the Moyal product and are thus denoted by the same symbol
$*_{\la}$. This could be a bit confusing because these new products
may become singular when $\la\rar 0$.

\vst

Fixing the set of generators $\phi_{\al}(x)$ of the ideal sheaf
$J_X$, one may equivalently define the above objects as follows
$$
\cJ_r = \{f\in\f_M[[\la]]\,\mid f=\sum_{\al=1}^n \pi^*(\phi_{\al})*_{\la}
g_{\al}\ \mbox{for some}\ g_{\al}\in \f_M[[\la]]\}
$$
$$
\cJ_l = \{f\in\f_M[[\la]]\,\mid f=\sum_{\al=1}^n g_{\al} *_{\la} \pi^*(\phi_{\al})
\ \mbox{for some}\ g_{\al}\in \f_M[[\la]]\},
$$
$$
\cN_r=\{ f\in \f_M[[\la]]\, \mid  f*_{\la}\pi^*(\phi_{\al})\subset
\cJ_r
\ \mbox{for all}\ \al=1,\ldots,k\},
$$
$$
\cN_l=\{ f\in \f_M[[\la]]\, \mid \pi^*(\phi_{\al})*_{\la} f \subset \cJ_l
\ \mbox{for all}\ \al=1,\ldots,k\},
$$
the equivalence (i.e.\ independence on the choice of generators)
being due to the associativity of $*_{\la}$ and the following elementary equality
$$
\pi^*\left(g(x) \phi_{\al}(x)\right)= \pi^*(g(x))*_{\la}
\pi^*(\phi_{\al}(x)) = \pi^*(\phi_{\al}(x))*_{\la} \pi^*(g(x)),
\ \ \forall g(x)\in \f_{\R^n}.
$$

\vst

\paragraph{\bf 2.1. Lemma.} {\em The sheaves of $*_{\la}$-algebras
$\cP_X$ and $\cQ_{X}$ are canonically isomorphic.}

\vst

\Proof This statement follows almost immediately from an elementary
observation that
$$
\pi^*(f)*_{\la} g=g*_{-\la}\pi^*(f)
$$
for any $f\in \f_{\R^n}$ and any $g\in \f_M$. $\Box$

\vst

The passage from $\cP_X$ to $\cQ_X$ (and vice versa) is essentially
equivalent to the transformation $\la\rar -\la$.

\vst

Global sections of the sheaf of $*_{\la}$-algebras $\cQ_X$ play the role of
admissible observables for the affine space $X$. The induced product $*_{\la}$
allows us to define, at least in principle, a spectral theory of observables via the
star-exponential \cite{F} and hence gives us the means to study
quantum mechanical models  on the background of a (singular) affine space
$X$. The physical interpretation of such models is very obscure --- it
is not even clear which observable would correspond to a ``free
particle" moving on $X$!

\vst

\paragraph{\bf 3. Quantization of the $n$-tuple point.}
The Moyal product on the cotangent bundle $M=\Omega^1 \R$ to the
real line is given explicitly by
$$
f*_{\la} g= \sum_{i=0}^{\infty} \frac{\la}{n!} \sum_{k=0}^i
\binom{i}{k} (-1)^k \frac{\p^i f}{\p x^{i-k}\p p^k} \frac{\p^i
g}{\p x^k \p p^{i-k}},
$$
where $f,g\in \f_M$.

\vst

Let $X$ be the $n$-tuple point in $\R$ given by the equation
$$
x^n=0.
$$
Then the left ideal is
$$
\cJ_l=\f_M *_{\la} x^n =\{f(x,p)*_{\la} x^n \ \mbox{for some smooth
function}\ f(x,p)\}
$$
and the associated normalizer is
$$
\cN_l = \{ f(x,p) \mid x^n*_{\la} f(x,p) = g(x,p)*_{\la} x^n
\ \mbox{for some smooth
function}\ g(x,p)\}.
$$

Any element $h(x,p)$ of the quotient $\cQ_X=\cN_l/\cJ_l$ can be uniquely
represented
as a sum
$$
h(x,p) = h_0 + h_1*_{\la}x + \ldots + h_{n-1}*_{\la} x^{n-1},
$$
for some functions of one variable $h_i=h_i(p)$, $i=0,1,\ldots,n-1$.
Moreover, this sum must satisfy the equation
$$
x^n*_{\la} h =0 \bmod \cJ_l,
$$
or equivalently
$$
\sum_{i=0}^{n-1}x^n *_{\la} h_i *_{\la} x^i =0 \bmod \cJ_l.
$$

\vst

\paragraph{\bf 3.1. Lemma.} {\em For any smooth function $g=g(x,p)$,
$$
x^i *_{\la} g = g*_{\la} x^i + \sum_{k=1}^i \binom{i}{k} (2\la)^k
\frac{\p^k g}{\p p^k}*_{\la} x^{i-k}.
$$
}

\vst

\noindent{\bf Proof} is straightforward.

\vst

\paragraph{\bf 3.2. Proposition.}{\em The ring of observables $\cQ_X$
consists of all possible sums $\sum_{i=0}^n h_i
*_{\la} x^i$, where the smooth functions $h_i=h_i(p)$ are given
explicitly by
$$
h_i =\sum_{k=0}^{n-1} a_{i,k}p^k +\frac{1}{(2\la)^i}
\sum_{k=0}^{i-1}\frac{(2\la)^k}{(i-k)!}
\left[\sum_{j=0}^{i-k-1}(-1)^{j+1} \binom{i-k-1}{j} a_{k,n-i+k+j} p^{n+j}\right]
$$
with $a_{i,k}$, $0\leq i,k\leq n-1$, being arbitrary constants.}

\vst

\Proof By Lemma~3.1,
$$
x^n*_{\la} h_i = \sum_{k=1}^n \binom{n}{k} (2\la)^k
\frac{\p^k h_i}{\p p^k}*_{\la} x^{n-k} \, \bmod \cJ_l.
$$
Hence, for any $h=\sum_{i=0}^{n-1}h_i*_{\la} x^i \in \cQ_X$,
\Beqr
0\bmod \cJ_l &=& \sum_{i=0}^{n-1}x^n *_{\la} h_i *_{\la} x^i\\
&=& \sum_{i=0}^{n-1} \sum_{k=1}^n
\binom{n}{k} (2\la)^k\frac{\p^k h_i}{\p p^k}*_{\la} x^{n-k+i}\\
&=&\sum_{i=0}^{n-1} \sum_{k=i+1}^n
\binom{n}{k} (2\la)^k\frac{\p^k h_i}{\p p^k}*_{\la} x^{n-k+i}\\
&=&\sum_{i=0}^{n-1} \sum_{l=i}^{n-1}
\binom{n}{l-i} (2\la)^{n-l+i}k\frac{\p^{n-l+i} h_i}{\p p^{n-l+i}}*_{\la}
x^{l}\\
&=&\sum_{l=0}^{n-1}\left[ \sum_{i=0}^{l}
\binom{n}{l-i} (2\la)^{n-l+i}k\frac{\p^{n-l+i} h_i}{\p p^{n-l+i}}\right]*_{\la} x^{l}
\Eeqr
implying that the functions $h_i=h_i(p)$ are solutions of the
following system of differential equations
$$
\sum_{i=0}^l \binom{n}{i-l} (2\la)^i \frac{d^{n-l+i}h_i}{dx^{n-l+i}} =
0, \ \ \
\ \ \ l=0,1,\ldots,n-1.
$$
The general solution of this system can be found by induction
and is precisely the one given in Proposition~3.2. $\Box$

\vst

If we assume that the deformation parameter $\la$ takes values in a field $\K$, which
is either $\R$ or $\C$, then we have to work with functions on $M$ with values in $\K$
and hence to view the constants of integration $a_{i,k}$ in
Proposition~3.2 as elements of $\K$.

\vst

\paragraph{\bf 3.3. Corollary.} {\em If $X$ is the $n$-tuple
point, then $\dim_{\K}\cQ_X=n^2$}.

\vst

\paragraph{\bf 3.4. Example.} If $n=2$, then a typical element $h$
of $\cQ_X$ is of the form
$$
h=a+bp+ \left(c + dp - b\frac{p^2}{2\la}\right)*_{\la} x,
$$
where $a,b,c$ and $d$ are constants, and the induced Moyal product is given by
\Beqr
h*_{\la}\tilde{h} &=& [a\tilde{a} + 2\la c\tilde{b}] + [a\tilde{b}
+ b\tilde{a} + 2\la d\tilde{b}]p \\
&=& \left[(a\tilde{c}+ c\tilde{a} +2\la c\tilde{d}) +
(a\tilde{d} + d\tilde{a} + b\tilde{c}-c\tilde{b} +2\la d\tilde{d})p - (a\tilde{b}
+b\tilde{a} +2\la d\tilde{b})\frac{p^2}{2\la}\right]*_{\la} x.
\Eeqr
Defining the map
$$
\Ba{rccc}
\psi: & \cQ_X & \lon & \mbox{Mat}_{\K}(2,2)\\
&        h & \lon & \left(\Ba{cc} a +2\la d &  b\\ 2\la c & 2\la a \Ea
\right)
\Ea
$$
one gets
$$
\psi(h *_{\la} \tilde{h}) =\psi(h) \cdot \psi(\tilde{h})
$$
where $\cdot$ stands for the usual matrix multiplication.
Hence $\psi$ identifies $\cQ_X$ together with the induced Moyal product $*_{\la}$
with the algebra of $2\times 2$-matrices.

\vst

\paragraph{\bf 3.5. Theorem.} {\em  If $X$ is the $n$-tuple
point, then the quantum algebra $(\cQ_X,*_{\la})$ is canonically
isomorphic, for any $\la\neq 0$, to the algebra of $n\times
n$-matrices.}

\vst

\Proof Let us define the $n$-dimensional vector space
$$
V=\mbox{span}_{\K}\left(e_0=x^{n-1}, e_1=p*_{\la} x^{n-1}, \ldots,
e_{n-1}=p^{n-1}*_{\la} x^{n-1}\right).
$$
Since the induced Moyal product is associative, we have, for any
$h\in \cQ_X$ and any $k=0,1,\ldots,n-1$,
\Beqr
h*_{\la}(p^k *_{\la} x^{n-1}) &=&
\sum_{i=0}^{n-1}
h_i *_{\la}(x^i *_{\la} p^k)*_{\la} x^{n-1} \\
&=&    \sum_{i=0}^{n-1}\sum_{l=0}^i
(2\la)^l \binom{i}{l}(h_i *_{\la} \frac{d^l p^k}{d p^l})*_{\la} x^{n+i-l-1}       \\
&=_{\bmod{\cJ_l}}&  \sum_{i=0}^{k}
(2\la)^i (h_i \frac{d^i p^k}{d p^i})*_{\la} x^{n-1}\\
&=_{\bmod{\cJ_l}}& g_k(p)*_{\la} x^{n-1},
\Eeqr
where $g_k(p)=\sum_{i=0}^k (2\la)^i h_i \frac{d^i p^k}{dp^i}$. We
claim that each function $g_k(p)$, $k=0,1,\ldots,n-1$, is a polynomial in $p$ of order
at most $n-1$. Indeed,
\Beqr
\frac{d^n g_k(p)}{d p^n} &=& \sum_{i=0}^k \sum_{j=0}^n (2\la)^i \binom{n}{j}
\frac{d^{i+j} p^k}{dp^{i+j}} \frac{d^{n-j} h_i}{dp^{n-j}} \\
&=& \sum_{l=0}^k \frac{d^l p^k}{d p^l} \sum_{i=0}^l (2\la)^i
\binom{n}{l-i} \frac{d^{n-l+i} h_i}{d p^{n-l+i}} \\
&=& 0,
\Eeqr
where we used the differential equations for $h_i(p)$ obtained
in the proof of Proposition~3.2.

\vst

Thus the resulting equality
$$
h*_{\la} e_a = \sum_{b=0}^{n-1} A^h_{ab} e_b, \ \
a=0,1,\ldots,n-1,
$$
with $A^h_{ab}\in \K$, defines a homomorphism from the algebra $(\cQ_X, *_{\la})$
to the algebra of $n\times n$ matrices,
$$
\Ba{rccc}
\psi: & \cQ_X & \lon & \mbox{Mat}_{\K}(n,n)\\
&        h & \lon & A^h_{ab}.
\Ea
$$

This homomorphism is injective. Indeed, if
$$
h*_{\la} p^k *_{\la} x^{n-1} = 0,  \ \ \ \ \forall k\in
\{0,1,\ldots,n-1\},
$$
then
$$
\sum_{i=0}^k (2\la)^i h_i \frac{d^i p^k}{dp^i}=0, \ \ \ \ \forall k\in
\{0,1,\ldots,n-1\},
$$
implying $h_i=0$.

Finally, the dimension counting implies that the map $\psi$
is an isomorphism. $\Box$

\paragraph{\bf 4. Concluding remarks.} (i) Theorem~3.5 implies that
deformation quantization of the $k$th order infinitesimal neighbourhood
of the embedding $\R^n\hook \R^{n+1}$ is equivalent to
introducing matrix valued functions on the associated phase space
$\R^{2n}=\Omega^1\R^n$
with the quantum product being the tensor product of the $2n$-dimensional
Moyal product and the matrix multiplication. Therefore, it is
natural to expect the appearance of matrix algebras in any
quantum theory where the background ``space-time" is thickened
into one extra dimension (cf.\ \cite{B}).

\vst

(ii) The construction in subsection~3 can be easily generalized
from the affine analytic subspaces, $X\hook \R^n$,  to analytic
subspaces of arbitrary ambient manifolds, $X\hook Y$, provided the cotangent
bundles $\Omega^1 Y$ come equipped with the star products $*_{\la}$
satisfying the equality
$$
\pi^*\left(f(x) g(x)\right)= \pi^*(f(x))*_{\la}
\pi^*(g(x)) = \pi^*(g(x))*_{\la} \pi^*(f(x)),
\ \ \forall f(x), g(x)\in \f_{Y},
$$
where $\pi: \Omega^1 Y\rar Y$ is the natural projection. This
will be the case, for example, if $*_{\la}$ comes from the torsion-free affine
connection on $Y$ (lifted to $\Omega^1Y$) via the Fedosov
construction \cite{Fe}.

\vst

(iii) The computations of algebras of quantum observables for other examples
mentioned in the introduction, say for
the ``tick" $X_2=\{(x,y)\in \R^2 \mid y^2-x^3=0\}$ or the
real line with one double point $X_3=\{(x,y)\in \R^2 \mid xy=0,
y^{2}=0\}$, are much easier than the one we did for
the $n$-tuple point. We leave the details to the interested
reader.
\pagebreak

\vse

{\small

{\sc
\begin{tabular}{l}
Department of Mathematics\\
University of Glasgow\\
15 University Gardens \\
 Glasgow G12 8QW, UK
\end{tabular}
}

\end{document}